# Lacunary Sections for Locally Compact Groupoids

by

Arlan Ramsay

University of Colorado

*Abstract:* It is proved that every second countable locally Hausdorff and locally compact continuous groupoid has a Borel set of units that meets every orbit and is what is called "lacunary," a property that implies that the intersection with every orbit is countable.

**Introduction**

Peter Forrest [4] studied free actions of $\mathbb{R}^n$ and showed that in the presence of a quasiinvariant measure, there is a complete lacunary section, i.e., a Borel set $Y$ whose saturation has a complement of measure 0, such that there is a neighborhood $U$ of the origin in $\mathbb{R}^n$ for which $y \in Y$ implies that $(y \cdot U) \cap Y = \{y\}$. Then $Y$ never contains more than countably many points in a single orbit. Applying Forrest's result in the case $n = 1$ recovers the Ambrose-Kakutani Theorem that every ergodic flow is "built under a function over a transformation." Hence his result is a generalization of that theorem. The main technique in [4] has been used in several papers to prove similar theorems in various contexts [3,6,7].

The paper of Kechris [6] proves the existence of lacunary sets that meet *all* orbits, not just almost all orbits, and does this for Borel actions of arbitrary second countable locally compact groups. This strengthens Forrest's result in two ways. First of all the group action need not be free, and besides that his result holds in 'the Borel category', which is stronger than the earlier results in 'the measure category'. In particular, the Kechris result provides a further generalization of the Ambrose-Kakutani Theorem and strengthens it so that the countable section is avialable for all $G$-quasiinvariant measures simultaneously.

Kechris reduces the proof to the case of continuous actions on compact spaces by using the (compact) universal '$G$-space' of Varadarajan [11]. Thus the existence of complete lacunary sections for Borel actions of locally compact groups follows from their existence for locally compact groupoids (defined below), which will be proved here. The paper of Kechris motivated the present proof, and the proof of [6] does generalize to groupoids. However, the proof given here follows ideas introduced in [5] and extended to groupoids in [10]. The theorem will be stated and comments made on it after some necessary definitions.

Speaking only algebraically at first, to make a *groupoid $G$ on a set $X$* we begin with *range* and *source* maps $r$ and $s$ from $G$ onto $X$ in terms of which (composability or) multiplicability of groupoid elements is defined. The multiplication is only defined for certain pairs and is a function from $G^{(2)} = \{(\gamma, \gamma') : s(\gamma) = r(\gamma')\}$ to $G$. We write $\gamma\gamma'$ for the product when it is defined. It is required that $r(\gamma\gamma') = r(\gamma)$ and $s(\gamma\gamma') = s(\gamma')$. Then *associativity* is formulated as follows: if $s(\gamma) = r(\gamma')$ and $s(\gamma') = r(\gamma'')$, then both products $\gamma(\gamma'\gamma'')$ and $(\gamma\gamma')\gamma''$ are defined and we simply require them to be equal. We require *identities* by saying that for each $x \in X$ there is a unique element $i_x$ of $G$ such that $r(i_x) = s(i_x) = x$ and so that if $\gamma \in G$ then $\gamma i_{s(\gamma)} = i_{r(\gamma)}\gamma = \gamma$. We identify $x$ with $i_x$, and therefore regard $X$ as a subset of $G$. For each $\gamma \in G$, there must be a unique $\gamma' \in G$ such that $s(\gamma) = r(\gamma')$ and $r(\gamma) = s(\gamma')$, $\gamma\gamma' = r(\gamma)$ and $\gamma'\gamma = s(\gamma)$. This element is called the *inverse* of $\gamma$ and written $\gamma^{-1}$. The shortest way to define *groupoid* is to say it means a small category with inverses.

To get a groupoid from a transformation group, suppose that a group $H$ acts on the right on a set $X$. Set $G = \{(x, h, y) \in X \times H \times X : xh = y\}$, $r(x, h, y) = x$, $s(x, h, y) = y$, $(x, h, y)(y, h', z) = (x, hh', z)$, $(x, h, y)^{-1} = (y, h^{-1}, x)$.

Now suppose that $X$ is second countable, locally compact and Hausdorff, that $G$ is second countable, locally compact and locally Hausdorff, and that $\{i_x : x \in X\}$ is closed in $G$ and homeomorphic to $X$, so that the identification is topological as well as set theoretic. We will call $G$ *continuous* if all the mappings involved in the definition are continuous. Cf. [7,9]. We call $G$ *topological* if it is continuous and the mapping $r$ is open. We will work with continuous groupoids. Second countability combined with local compactness implies that the groupoid is a countable union of compact sets, i.e., is $\sigma$-compact.

We write $\theta$ for the mapping $(r, s)$ from $G$ into $X \times X$. The range of $\theta$ is the equivalence relation $R_G$ induced by $G$ on $X$. The groupoid $G$ is called *principal* provided that $\theta$ is one-one, so that $G$ is isomorphic to $R_G$ as a groupoid. Since $G$ is $\sigma$-compact, so is $R_G$ and hence it is a Borel set in $X \times X$. For $x \in X$, the set $[x] = s(xG)$ is called the *orbit* of $x$. The orbits partition $X$. We write $G(x)$ for $\{\gamma \in G : s(\gamma) = r(\gamma) = x\}$ and $G'$ for $\{\gamma \in G : s(\gamma) = r(\gamma)\}$. The first is a group, called the stabilizer or isotropy of $x$. The second is a closed set called the *isotropy bundle,* and it will be proved to be a continuous groupoid. The function $S$ assigning to $x \in X$ the group $G(x)$ is called the *isotropy map,* and will play a significant role in the proof.

Since $r$ may not be open, we cannot prove that $AB$, defined to be the image of $(A \times B) \cap G^{(2)}$ under multiplication, is open whenever $A$ and $B$ are open. However, it is true that $AB$ is closed if one of the sets is closed and the other is compact. If $A$ is a subset of $X$ and $B$ is a subset of $G$, then we write $[A]B$ for the set $s(AB)$ in $X$. Thus for a point $x \in X$, we have $[x] = [x]G$.

**Definition.** ([3,4,6]) A Borel set $Y \subseteq X$ is a (partial) *countable section* if for every $x \in X$ the orbit $[x]$ meets $Y$ in a countable set. A countable section $Y$ is *complete* if $Y$ contains a point of every orbit. The Borel set $Y$ is a *complete lacunary section* iff it contains a point of every orbit and there is a set $U$ that is open in $G$ and contains $X$, such that if $y \in Y$ then $([y]U) \cap Y = \{y\}$.

We will see below that a lacunary section is always a countable section. The main theorem of this paper is as follows:

**Theorem A.** If $G$ is a continuous locally Hausdorff locally compact groupoid on $X$, then $X$ contains a complete lacunary section $Y$.

The most familiar non-Hausdorff groupoid is the holonomy groupoid of a foliation, where the existence of lacunary sections comes directly from the definition of foliation. The present theorem may be regarded as proof that local Hausdorffness is sufficient in general.

Even if we started with a topological groupoid, the proof would involve groupoids in which the mapping $r$ might not be open. Moreover, there are potential applications of this form of the theorem, and it does allow the section to be placed in any closed set that meets every orbit. Some applications are expected to be similar to applications of the Ambrose-Kakutani type theorems for transformation groups. For example, one may be able to understand amenability for general locally compact groupoids by reducing some questions to the case of countable orbits.

Further, one may speculate about uses of groupoids that have no (continuous) Haar system. They do have unitary representations, making for possible applications. The existence of lacunary sections may well be a useful tool in working with such groupoids.

**Acknowledgements:** The author did much of the work on this paper while enjoying the kind hospitality of the School of Mathematics at Flinders University. He thanks Professor William Moran for the opportunity to pursue research there. He also thanks the University of Colorado for the Sabbatical which provided the freedom to travel.

## 1. Proof of Theorem A

The plan of the proof is to use a modification of Forrest's method obtained by using some ideas of Gootman and Rosenberg [5] as adapted to groupoids by Renault [10]. The idea is that a subset of $G$ that contains the set of units in its interior can induce an equivalence relation on appropriate small neighborhoods.

If $G$ were determined by a transformation group, we could choose a compact neighborhood of the identity in the group to use in making standardized neighborhoods in the orbits. Lacking that, we proceed by choosing a symmetric conditionally compact set $D$ in $G$ whose interior contains $X$. Conditional compactness means that if $B$ is a compact set in $X$, then the sets $BD$ and $DB$ are compact. (It follows that $BDB$ is compact because the closedness of $B$ implies that $BDB$ is a closed set in $BD$.) Such a set $D$ can always be found because $G$ and $X$ are both locally compact and second countable. We will find a Borel section whose lacunarity is established by $D^\circ$, the interior of $D$.

Set $R_D = \theta(D)$, a subset of the relation $R_G$. The strategy is to seek subsets of $X$ on which the relation $R_D$ provides an equivalence relation that has a transversal: a Borel set that meets each equivalence class exactly once. If there are enough such sets and they behave well, the proof can be completed. Since $D$ is symmetric and contains $X$, the relation $R_D$ is symmetric and reflexive. Transitivity is what we seek, and it can be achieved by using the ideas of Gootman-Rosenberg and of Renault [5,10]. If $x$ is a point of continuity for the isotropy map, we can show that $x$ has a neighborhood on which $R_D$ is an equivalence relation, and compact subsets of such a neighborhood will contain Borel transversals.

At least in the Hausdorff case, the result we want has been proved for principal groupoids, so the presence of non-trivial isotropy provides the main barrier we must surmount. Thus we begin by showing that the isotropy bundle is Hausdorff and locally compact even though $G$ is not required to be Hausdorff.

**Lemma 1.** If $G$ is a locally Hausdorff locally compact continuous groupoid, then the isotropy bundle, $G'$, is closed in $G$ and hence locally compact. In addition, $G'$ is Hausdorff and is a continuous groupoid.

**Proof.** The fact that $G'$ is closed in $G$ follows from the fact that $X$ is Hausdorff and $\theta$ is continuous from $G$ into $X \times X$. Continuity of the operations is clear. The fact that $G$ is locally Hausdorff implies that $G'$ is locally Hausdorff and the fact that $G'$ is closed in $G$ implies that $G'$ is locally compact. For each $x \in X$ it follows that that $xG'$ is also locally Hausdorff and locally compact. Furthermore, $xG'$ is a group and the group operations are continuous. Hence $xG'$ is a Hausdorff locally compact group. To prove that $G'$ is Hausdorff, we begin with points $\gamma$ and $\gamma'$ for which $r(\gamma) \neq r(\gamma')$. Disjoint neighborhoods of $r(\gamma)$ and $r(\gamma')$ have inverse images under $r$ that are disjoint neighborhoods of $\gamma$ and $\gamma'$. Now suppose that $r(\gamma) = r(\gamma') = x$ but $\gamma \neq \gamma'$. Then there are compact Hausdorff sets $K$ and $K'$ that are neighborhoods of $\gamma$ and $\gamma'$ respectively and such that $K \cap xG'$ and $K' \cap xG'$ are disjoint. If $K$ and $K'$ are disjoint, the proof is complete. If $K \cap K' \neq \emptyset$, then $r(K \cap K')$ is a compact set in $X$. Since $X$ is Hausdorff, $r(K \cap K')$ is closed. Hence $r^{-1}(r(K \cap K'))$ is a closed set in $G'$ that is disjoint from $xG'$. Removing $r^{-1}(r(K \cap K'))$ from both $K$ and $K'$ provides disjoint neighborhoods of $\gamma$ and $\gamma'$ in $G'$.

In order to use the ideas of [5,10], we need to study the isotropy map, $S$, that assigns to each $x \in X$ the isotropy group $G(x)$. The values of $S$ are locally compact groups, but it is equally important that they are closed sets in $G'$. We write $\mathcal{C}(G')$ for the set of closed subsets of $G'$, equipped with the Fell topology. A basis for this topology is made of sets defined as follows: let $K$ be a compact set, let $U_1, U_2, \ldots, U_n$ be open sets, and set $\mathcal{U}(K, U_1, U_2, \ldots, U_n) = \{F \in \mathcal{C}(G') : F \cap K = \emptyset \text{ and } F \cap U_i \neq \emptyset$ for $i = 1, \ldots, n\}$. Using the fact that $G'$ is Hausdorff, locally compact and second countable, we know from [1] that $\mathcal{C}(G')$ is metrizable and that the topology is the weak topology determined by a sequence of real valued functions defined as follows: Let $\gamma_1, \gamma_2, \ldots$ be a sequence that is dense in $G'$. Let $d$ be a metric that determines the topology on $G'$. For each $i$, choose a positive number $a_i$ such that the closed ball of radius $2a_i$ centered at $\gamma_i$ is compact. Then define $d_i(B) = \min(a_i, d(\gamma_i, B))$, where the second component is the distance from $\gamma_i$ to $B$. We will show that in our setting such a function composed with the isotropy map gives a semicontinuous function. This fact implies that $S$ always has a dense $G_\delta$ of continuity points.

For each $x \in X$, $S(x) = xG'$, so the graph of $S$ is the set $\{(x, \gamma) \in G^{(2)} : s(\gamma) = x\}$, which is closed in $X \times G'$ because $G^{(2)}$ is closed in $G \times G$ and $s$ is continuous. The fact that the graph is closed suffices to prove the lower semicontinuity of $d_i(S(x))$ as a function of $x$, which implies what we needed. Then properties of real valued semicontinuous functions give us the next lemma.

**Lemma 2.** If $Z$ is a closed subset of $X$, the restriction of the isotropy map to $Z$ has a dense $G_\delta$ set of continuity points.

The next three lemmas are found in [10]. We give the proofs because we are not requiring $r$ to be open, as Renault did.

**Lemma 3.** Let $G$ be a continuous groupoid and let $x \in X$. The following properties are equivalent.

i) The isotropy map is continuous at $x$.

ii) The restriction of $r$ to $G'$ is open at every $\gamma \in G(x)$ (i.e., sends every neighborhood

of $\gamma$ onto a neighborhood of $x$.)

iii) For every $\gamma \in G(x)$ and every neighborhood $U$ of $\gamma$ in $G$, $G'U$ is a neighborhood of $G(x)$ in $G$.

**Proof.** If $K$ is compact in $G$ and $U_1, U_2, \ldots, U_n$ are open sets in $G$, the inverse image under $S$ of the basic open set $\mathcal{U}(K, U_1, U_2, \ldots, U_n)$ is

$$r(K \cap G')^c \cap r(U_1 \cap G') \cap \cdots \cap r(U_n \cap G').$$

Thus continuity of $S$ at $x$ implies that if $\gamma \in G(x)$ and $U$ is a neighborhood of $\gamma$ in $G$, then $r(U \cap G')$ is a neighborhood of $x$, i.e., $r|G'$ is open at each $\gamma \in G(x)$. Thus i) implies ii). We know that $r(K \cap G')$ is compact and hence closed in $X$ whenever $K$ is compact in $G$, for any continuous groupoid. Thus the formula above makes it clear that ii) implies i). To see that iii) implies ii), we just notice that $r(U \cap G') = G'U \cap X$. To prove that ii) implies iii), we form the semidirect product $G' \ltimes G = \{(\gamma', \gamma) \in G' \times G \colon s(\gamma') = r(\gamma)\}$. This is a groupoid if we let $G'$ act on $G$. Let $\pi_2$ denote the projection of $G' \ltimes G$ onto the second factor and define $\pi_1(\gamma', \gamma) = \gamma'^{-1}\gamma$, obtaining a continuous function from $G' \ltimes G$ to $G$. In preparation for the proof that ii) implies iii) we want to see that $\pi_2$ is open at every point of $G(x) \ltimes xG$, under the assumption that $r|G'$ is open at every point of $G(x)$. Let $\gamma' \in G(x)$, let $\gamma \in xG$ and let $U$ be a neighborhood of $(\gamma', \gamma)$ in $G' \ltimes G$. Then there are open sets $V$ and $W$ in $G$ containing $\gamma'$ and $\gamma$ respectively, so that $(V \cap G') \ltimes W \subseteq U$. Hence $\pi_2(U)$ contains $W \cap (r^{-1}(r(V \cap G')))$. Since $r|G'$ is open at $\gamma'$, we see that $\pi_2(U)$ is a neighborhood of $\gamma$ in $G$. Now, if $\gamma$ is a point of $G(x)$ and $U$ is a neighborhood of $\gamma$ in $G$, then $\pi_1^{-1}(U)$ is a neighborhood of every point of $G' \ltimes G$ that $\pi_1$ carries to $\gamma$, namely every point of $\{(\gamma', \gamma'\gamma) : \gamma' \in G(x)\}$. We see that $G'U = \pi_2(\pi_1^{-1}(U))$, and the fact that $\pi_2$ is open at every point of $G(x) \ltimes xG$ completes the proof.

For $K \subseteq G$, we write $R_K$ for the set $\theta(K)$, a relation on $X$. We will see that the restriction of $R_K$ to suitable subsets is an equivalence relation if $K$ is conditionally compact, symmetric and contains $X$ in its interior. Whenever $K$ is conditionally compact the relation is closed. (Let $(\gamma_i)$ be a net in $K$ such that $\theta(\gamma_i)$ converges to a point $(y, z)$ in $X \times X$. By passing to a subnet, we may suppose that $(\gamma_i)$ converges to a point of $K$, and it follows that $(y, z) \in R_K$.)

**Lemma 4.** Let $K$ be a conditionally compact subset of $G$, $x \in X$ and $N$ an open neighborhood of $xKx$ in $G$. Then there is a neighborhood $V$ of $x$ in $X$ such that $VKV \subseteq N$.

**Proof.** If not, for every compact neighborhood $V$ of $x$ in $X$ the set $VKV \cap N^c$ is not empty, so we can choose a point $\gamma_V$ in that set. For any particular compact neighborhood $W$ of $x$, there is a subnet of the $\gamma_V$'s that converges in $WKW$. If that limit point is $\beta_W$, we know that $r(\beta_W) = s(\beta_W) = x$, and $\beta_W \in WKW \cap N^c$. But then, $\beta_W \in xKx \cap N^c$, which is impossible.

**Lemma 5.** (Cf. [5, Theorem 1.4]) Let $x$ be a point of continuity of the isotropy and let $K$ be a conditionally compact symmetric neighborhood of $x$ in $G$. Then there is an open neighborhood $V$ of $x$ in $X$ such that $R_K$ restricted to $V$ is a closed equivalence relation.

**Proof.** It follows from Lemma 3 that $G'K$ is a neighborhood of $G(x)$ in $G$. Since $xKx$ is compact and a product of two of its elements is always in $G(x)$, there is an open neighborhood $N$ of $xKx$ such that $N^2 \subseteq G'K$. Then there is an open neighborhood $V$ of $x$ in $X$ such that $VKV \subseteq N$ and $V \subseteq N$. Since $N^2 \subseteq G'K$, we have

$$R_K \cap (V \times V) \subseteq R_N \cap (V \times V) \subseteq R_{N^2} \cap (V \times V) \subseteq R_K \cap (V \times V).$$

The fact that $V \subseteq N$ implies that the relation $F_V := R_K \cap (V \times V)$ is reflexive. The symmetry of $K$ implies that $F_V$ is symmetric. It is transitive because of the containments written above. It is closed in $V \times V$ because $K$ is conditionally compact.

Of course Lemma 5 can be applied to the set $K = D$. The next lemma establishes some properties of the equivalence relation $F_V$ defined in proof of Lemma 5.

**Lemma 6.** Let $x$ be a point of continuity of the isotropy map, and let $V$ be an open neighborhood of $x$ such that $R_D \cap (V \times V)$ is a closed equivalence relation. Let $B$ be a compact subset of $V$ and set $F_B = R_D|B$.

i) $F_B$ has a Borel transversal, $T$.

ii) If $T$ is a subset of $B$ that contains at most one point of each $F_B$ class, then $T$ is $D°$-lacunary, so that $T \cap [x]$ is a discrete subset of $[x]$ when $[x]$ is given the quotient topology of $xG$ relative to the map $s$. Hence $T \cap [x]$ is countable.

**Proof.** The space $B$ is a compact metric space, and $F_B$ is a closed equivalence relation, so the quotient space is also a compact metric space. Also, the saturation of every closed set is closed. Hence either the Federer-Morse transversal lemma or the Dixmier transversal lemma can be used to prove part i) [2]. For part ii), we use the fact that $G(x)$ is a locally compact group acting topologically on $xG$, so that the saturation of an open set is open. We also know that $s$ is continuous from $xG$ to $[x]$, so the image of a compact set is compact. Now, suppose that $y \in [x] \cap T$. Then $yD = r^{-1}(y) \cap D$ is a compact neighborhood of $y$ in $yG$, so $[y]D = s(yD)$ is a compact neighborhood of $y$ in $[y] = [x]$. The quotient topology inherited from $yG$ is the same as that inherited from $xG$. If $z \neq y$ is another element of $[x] \cap T$, then $[y]D \cap [z]D = \emptyset$ because $F_B$ is an equivalence relation. Hence each element of $T \cap [x]$ has a neighborhood in the quotient topology on $[x]$ that contains no other element of $T \cap [x]$, as needed.

**Final steps in the proof:** Begin with a continuous groupoid $G$ on $X$. Let $S$ be the isotropy map. Take $D$ to be a symmetric conditionally compact set that contains $X$ in its interior. We will exhaust $X$ using a transfinite induction, and then prove that countably many sets suffice by using a countable basis for the topology of $X$.

Consider the collection $\mathcal{W}$ of open sets $W$ in $X$, such that every point $x \in W$ has an open neighborhood $V$ for which $R_D \cap (V \times V)$ is a closed equivalence relation. It is clear that $\mathcal{W}$ is closed under arbitrary unions, so it contains a largest element, $W_0$. By Lemmas 2 and 4, we know that $W_0$ is dense in $X$. Then we set $X_1 = X \backslash W_0$, define $G_1 = G|X_1$, $S_1 = S|X_1$ and observe that $S_1$ has a dense set of continuity points. Hence, as before, there is a set, $W_1$, dense in $X_1$, and relatively open, such that each point of $W_1$ has an open neighborhood $V$ in $X_1$ for which $R_D \cap (V \times V)$ is a closed equivalence relation.

For the induction step, we suppose that we have constructed closed sets $X_\beta$ for $\beta < \alpha$, so that $\beta < \gamma$ implies $X_\beta \supset X_\gamma$. If $\alpha$ is a limit ordinal, set $X_\alpha = \bigcap\{X_\beta : \beta < \alpha\}$. If $\alpha$ is the successor of $\beta$, find an open dense set $W_\beta$ in $X_\beta$ such that each point of $W_\beta$ has a neighborhood $V$ in $X_\beta$ for which $R_D \cap (V \times V)$ is a closed equivalence relation. Then take $X_\alpha = X_\beta \backslash W_\beta$.

To prove that this induction process reaches $X_\alpha = \emptyset$ for some countable ordinal, we first remark that $X_\alpha = \emptyset$ for $\alpha$ greater than the first ordinal whose cardinality is that of $X$. Next, let $\mathcal{U}$ be a countable basis for the topology of $X$. We denote by $\mathcal{U}_\alpha$ the subset of $\mathcal{U}$ consisting of those basis elements that are disjoint from $X_\alpha$ so that $X_\alpha = \emptyset$ if $\mathcal{U}_\alpha = \mathcal{U}$. Notice that $\mathcal{U}_\alpha \subset \mathcal{U}_\beta$ whenever $X_\alpha \supset X_\beta$, and there can only be countably many strict increases. This establishes the desired countability.

For each $\alpha$, the set $W_\alpha$ is a countable union of sets open in $X_\alpha$ on which $R_D$ is a closed equivalence relation, and each one of those smaller open sets is a countable union of compact sets. Hence, there is a sequence $B_1, B_2, \ldots$ of compact sets whose union is $X$ such that for each $n$ $R_D \cap (B_n \times B_n)$ is a closed equivalence relation. For each $n$, let $T_n$ be a Borel transversal of $R_D \cap (B_n \times B_n)$. Since $B_n G$ is closed in $G$, it is a countable union of compact sets, and hence the same is true for $[B_n] = s(B_n G)$, the saturation of $B_n$ in $X$. Define

$$Y = \bigcup_{n \geq 1} T_n \backslash \left( \bigcup_{k < n} [B_n] \right).$$

This set does what is claimed in Theorem A.